\documentclass[12pt]{amsart}
\usepackage{amscd,amsmath,amssymb,amsfonts}
\usepackage[cmtip, all]{xy}
\theoremstyle{plain}
\newtheorem{thm}{Theorem}

\theoremstyle{definition}

\numberwithin{thm}{section}
\numberwithin{equation}{section}

\newcommand{\ml}[2]{\begin{multline}\label{#1}#2 \end{multline}}
\newcommand{\ga}[2]{\begin{gather}\label{#1}#2 \end{gather}}

\newcommand{\surj}{\twoheadrightarrow}

\newcommand{\sA}{{\mathcal A}}

\newcommand{\sF}{{\mathcal F}}

\newcommand{\C}{{\mathbb C}}

\newcommand{\F}{{\mathbb F}}

\renewcommand{\P}{{\mathbb P}}
\newcommand{\Q}{{\mathbb Q}}

\begin{document}

\title[Frobenius Eigenvalues]{Eigenvalues of Frobenius acting on 
the $\ell$-adic cohomology of complete intersections of low 
degree}
\author{H\'el\`ene Esnault}
\address{
Universit\"at Essen, FB6, Mathematik, 45117 Essen, Germany}
\email{esnault@uni-essen.de}

\date{April 18, 2003}
\begin{abstract}
We show that the eigenvalues of Frobenius acting on $\ell$-adic cohomology 
of a complete intersection of low degree defined over the finite field $\F_q$
modulo the cohomology of the projective space
are divisible as algebraic integers by $q^\kappa$, where 
the natural number $\kappa$ is predicted by the theorem of Ax and Katz (\cite{Ka}) on the congruence for the number of rational points. \\ \ \\
{\bf Valeurs propres de Frobenius agissant sur la cohomologie $\ell$-adique
d'intersections compl\`etes de bas degr\'e.}\\ \ \\
{\bf R\'esum\'e}: Nous montrons que les valeurs propres de Frobenius agissant
sur la cohomologie $\ell$-adique d'une intersection compl\`ete de bas degr\'e
d\'efinie sur le corps fini  $\F_q$ modulo la cohomologie de
$\P^n$ sont divisibles en tant
qu'entiers alg\'ebriques  par $q^\kappa$, o\`u 
$\kappa$ est l'entier naturel pr\'edit par le th\'eor\`eme de Ax et
Katz (\cite{Ka}) sur la congruence pour le nombre de points rationnels. 

\end{abstract}
%\subjclass{Primary Algebraic Geometry}
\maketitle
\begin{quote}

\end{quote}
{\bf Version fran\c{c}aise abr\'eg\'ee}. 
Si $X\subset \P^n$ est une intersection 
compl\`ete de degr\'es $d_1\ge \ldots \ge d_r$ d\'efinie sur le corps fini
$\F_q$, et si $\kappa$ est la partie enti\`ere du nombre rationnel 
$\frac{n-d_2-\ldots -d_r}{d_1}$, nous montrons que les valeurs propres 
de Frobenius agissant sur la cohomologie $\ell$-adique 
relative 
$H^i(\overline{\P^n}, \overline{X}, \Q_\ell)$ 
sont divisibles par $q^\kappa$ en tant qu'entiers alg\'ebriques 
(Theorem \ref{mainthm}). Ici, $
\overline{\P^n}=\P^n\times_{\F_q} \overline{\F_q}, \
\overline{X}=X\times_{\F_q} \overline{\F_q}$
et $\overline{\F_q}$ est la cl\^oture alg\'ebrique de $\F_q$. 
Si $X$ est lisse, utilisant le th\'eor\`eme de 
Ax et Katz (\cite{Ka}), c'est une cons\'equence imm\'ediate
de la formule des traces de Grothendieck-Lefschetz 
\cite{Gr} et de ce que la cohomologie primitive de $\overline{X}$ est concentr\'ee en dimension moiti\'e \cite{DeSGA}. Si $X$ est singuli\`ere, cela n'est plus le cas. Si $\F_q$ est remplac\'e par un corps $k$ de caract\'eristique 0, 
l'assertion qui correspond dans la philosophie motivique 
au Theorem \ref{mainthm} est que le type de Hodge est $\ge \kappa$ sous les
m\^emes hypoth\`eses de degr\'e et d'intersection compl\`ete,  et est prouv\'ee dans  \cite{DD} pour les hypersurfaces et \cite{E} en g\'en\'eral. 
Cela dit, le type de Hodge est $\ge \kappa$ m\^eme sans hypoth\`ese d'intersection compl\`ete (\cite{ENS}), mais la preuve de cette note ne permet pas 
d'\^oter l'hypoth\`ese d'intersection compl\`ete dans le
Th\'eor\`eme \ref{mainthm}.

\section{Introduction}
If $X\subset \P^n$ is a projective variety 
defined over the finite field $k=\F_q$
by equations of degrees 
$d_1\ge d_2\ge \ldots \ge d_r$, then the fundamental theorem of Ax and Katz
\cite{Ka} asserts that the number of rational points of $X$ fulfills
\ga{1}{|X(\F_q)| \equiv |\P^n(\F_q)| \ \text{mod} \ q^\kappa.}
If $X$ is smooth, this is equivalent to saying that the eigenvalues of 
the Frobenius action on 
$H^i(\overline{\P^n}, \overline{X}, \Q_\ell)$ for all $i$ are of the shape
\ga{2}{q^\kappa\cdot {\rm algebraic \ integer}.}
Here $\overline{X}=X\times_{\F_q} \overline{\F_q}$  and $\overline{\F_q}$ is the algebraic closure of $\F_q$, and $\kappa$ is the integral part of $\frac{n-d_1-\ldots -d_r}{d_1}$. 
 Let us  briefly recall  why. 

Let $U=\P^n\setminus X$, and let $\zeta(U,t)$ be the zeta function of $U$
defined by its logarithmic derivative
\ga{3}{ 
\frac{\zeta'(U, t)}{\zeta(U,t)}=\sum_{s\ge 1} |U(\F_{q^s})|t^{s-1}.}
By the theorem of Dwork \cite{Dw}, we know that 
$\zeta(U, t)$  is a rational function 
\ga{4}{\zeta(X,t)\in \Q(t),}
while Grothendieck-Lefschetz trace formula \cite{Gr} gives a cohomological
formula
\ga{5}{\zeta(U,t)=\prod_{i=0}^{2\ \rm dim(U)} {\rm det}(1-F_it)^{(-1)^{i+1}},}
where $F_i$ is the arithmetic Frobenius acting on the compactly supported
$\ell$-adic cohomology  $H^i_c(\overline{U}, \Q_\ell)=
H^i(\overline{\P^n}, \overline{X}, \Q_\ell).$   
For $X$ smooth, the Weil conjectures 
\cite{DeWeI} 
assert that the eigenvalues in any complex
embedding $\Q_\ell \subset \C$
of the arithmetic Frobenius  acting on $H^i(\overline{X}, \Q_\ell)$
have absolute values $q^{\frac{i}{2}}$.
If we moreover assume that 
$X$ is a complete intersection, then we know
that only $H^{n-r}_c(\overline{X}, \Q_\ell)$ is 
not equal to $H^{n-r}(\overline{\P^n}, \Q_\ell)$
(\cite{DeInt})
and we don't need the Weil conjectures. 
The conclusion is that there is no possible 
cancellation between the numerator and the denominator of \eqref{5}. Thus
\eqref{1} for all finite field extensions $\F_{q^s}\supset \F_q$
is equivalent to \eqref{2}.

On the other hand, if we replace $k$ by a field of characteristic 0, keeping the same degree
assumption, then we know that the Hodge type of $X$, that is the largest integer $\mu$ such that the Hodge filtration in de Rham cohomology satisfies
\ga{6}{F^\mu H^i_{c, DR}(U)=H^i_{c, DR}(U)}
is $\ge \kappa$ (see \cite{DeSGA} in the smooth case, \cite{DD} for hypersurfaces, 
\cite{E} for complete intersections and \cite{ENS} for the general case). 
The purpose of this note is to show
\begin{thm} \label{mainthm}
Let $X\subset \P^n$ be a complete intersection defined over the finite field 
$k=\F_q$ by equations of degrees 
$d_1\ge d_2\ge \ldots \ge d_r$. Then the eigenvalues of 
the arithmetic Frobenius acting on $\ell$-adic cohomology $H^i(
\overline{\P^n}, \overline{X}, \Q_\ell)$ are
 of the shape $q^\kappa\cdot {\rm algebraic \ integer}$  for all $i$, 
where $\kappa $ is the integral part 
of $\frac{n-d_2-\ldots -d_r}{d_1}$.  
\end{thm}
The proof reduces the assertion to Deligne's integrality statement 
\cite{DeInt} which is 
true for all projective varieties, 
without supplementary assumptions. 
However, the reduction depends havily on the complete intersection assumption.

\noindent {\it Acknowledgements}. 
I thank Spencer Bloch, Pierre Deligne, Nick Katz,  Madhav Nori
 and Daqing Wan for discussions on 
this topic. In particular, Spencer Bloch warned that the 
Gysin theorem \ref{thm:gysin} needed in the proof
has to be justified, and Pierre Deligne provided a proof of it.  
In our induction, we originally reduced the theorem to 
\cite{DeInt}, Corollaire 5.5.3, (ii). Nick Katz observed that without 
changing the induction, it is enough to reduce to \cite{DeInt}, 
Corollaire 5.5.3, which is an easier statement.

\section{The proof of Theorem \ref{mainthm}}
 If $X$ is a complete intersection of dimension $n-r$,
Artin's vanishing theorem \cite{Ar}, Th\'eor\`eme 3.1 implies
\ga{2.1}{ H^i(\overline{X}, \Q_\ell)/H^i(\overline{\P^n}, \Q_\ell)=
H^{i+1}_c(\overline{U}, \Q_\ell)= 0 \ \text{for} \  i< {\rm dim} X=n-r.}  
By the Lefschetz trace formula \eqref{5}, we see, as already observed by D. Wan (see 
introduction of \cite{BEL}), that the theorem of Ax-Katz  \eqref{1} together with the divisibility assertion for $F_i$ for all $i$ but one, implies the divisibility assertion
for the last $F_i$. Thus we just have to show divisibility for $i>n-r$, or equivalently divisibility of  the $F_i$ acting on $H^i_c(\overline{U}, \Q_\ell)=H^i(\overline{\P^n}, 
\overline{X}, \Q_\ell)$ 
 for $i> n-r +1$. 
Let $A\cong \P^{n-1}$ be a linear hyperplane in general position. 
One has an exact sequence of 
\ml{2.2}{\ldots \to H^i_A(\overline{\P^n}, \overline{X}, \Q_\ell) \to 
H^i(\overline{\P^n}, 
\overline{X}, \Q_\ell) \\ \to 
H^i(\overline{\P^n\setminus A}, \overline{X\setminus X\cap A}, \Q_\ell) \to \ldots } 
which is compatible with the Frobenius action. 
For $i>n-r+1$, Artin's vanishing theorem \cite{Ar} again implies 
$H^i(\overline{\P^n\setminus A}, \overline{X\setminus X\cap A}, \Q_\ell)=0.$ Thus 
we are reduced to proving the assertion for the Frobenius action on 
$H^i_A(\overline{\P^n},\overline{X}, \Q_\ell)$ for $i>n-r+1$. 
One has
\begin{thm}[P. Deligne] \label{thm:gysin}
The Gysin homomorphism
\ga{}{H^{i-2}(\overline{A}, \overline{X\cap A}, \Q_\ell)(-1)\xrightarrow{{\rm Gysin}}
H^i_A(\overline{\P^n}, \overline{X}, \Q_\ell) \notag }
is an isomorphism of Frobenius modules 
 for $A$ in a non-trivial  open subset of the dual projective space 
$(\P^n)^\vee$. More generally, if $\sF$ is a $\ell$-adic sheaf 
on $\P^n$, then for $A$ in a 
non-trivial  open subset of the dual projective space 
$(\P^n)^\vee$, the Gysin homomorphism 
\ga{}{H^{i-2}(A, i^*\sF)(-1)\to H^i_A(\P^n, \sF),\notag}
where $i: A\to \P^n$ is the closed embedding, is an isomorphism. 
\end{thm}
\begin{proof}(P.  Deligne) We consider the universal family
\ga{2.3}{ \iota: \sA=\{(x, A)\in \P^n \times_k (\P^n)^\vee, x\in A \}\hookrightarrow 
\P^n \times_k (\P^n)^\vee}
together with the projections
\ga{2.4}{
 {\rm pr}_2:
\P^n \times_k (\P^n)^\vee \to (\P^n)^\vee, \
{\rm pr}_1:
\P^n \times_k (\P^n)^\vee \to (\P^n).
}
Since $\sF$ comes from $\P^n$, and 
$\sA\subset  \P^n\times_k (\P^n)^\vee$ is a smooth hypersurface such that 
$\sA\to \P^n$ is smooth, 
the classical Gysin homomorphism is an isomorphism
\ga{2.5}{\iota^* {\rm pr}_1^*\sF [-2](-1) \xrightarrow{ \cong} 
\iota^{!} {\rm pr}_1^*\sF.
}
By \cite{DeFin}, Corollaire 2.9 applied to ${\rm pr}_2$,
 there is a non-empty open  subset 
$V\subset (\P^n)^\vee$ such that 
for all $A\in V$ one has the base changes
\ga{2.6}{ i_A^*(\iota^* {\rm pr}_1^*\sF [-2](-1))=i^*\sF[-2](-1),
  \ i_A^*(\iota^{!}{\rm pr}_1^*\sF )= i^!\sF,
}
where $i_A: \sA\cap q^{-1} A=A\times \{A\} \hookrightarrow \P^n\times \{A\}$. 
Thus \eqref{2.5} and \eqref{2.6} imply that for such a $A$, 
$i^*\sF[-2](-1) \xrightarrow{\cong} i^!\sF$. This concludes the proof. 
\end{proof}
Thus for $i>n-r+1$, \eqref{2.2} becomes
\ga{2.7}{H^{i-2}(\overline{A}, \overline{X\cap A}, \Q_\ell)(-1) \surj 
H^i(\overline{\P^n}, 
\overline{X}, \Q_\ell)}
as a surjection of Frobenius modules. On the other hand, $X\cap A$ is 
a complete intersection of the same degrees $d_1,\ldots, d_r$ and we have
\ga{2.8}{\kappa(A\cap X)\le \kappa -1.}
Arguing by induction on $\kappa$, we are reduced to showing 
Theorem \ref{mainthm} 
for $\kappa=0$ and all $i$ (here notice in \eqref{2.7} the cohomological degree
drops by 2 in the Gysin map, thus we have to take 
into account the middle cohomology as well). 
Thus we have the assumptions of Theorem \ref{mainthm} for $\kappa \ge 0$. 
Then Theorem \ref{mainthm} simply says that the eigenvalues of $F_i$
are algebraic integrers, which is 
 Deligne's integrality theorem \cite{DeInt}, Corollaire 5.5.3.

\bibliographystyle{plain}

\renewcommand\refname{References}

\end{document}